\documentclass[10pt]{article}

\usepackage{authblk}
\usepackage{fullpage}

\usepackage{mathtools,amssymb,amsthm,amsfonts,wasysym}
\usepackage[maxbibnames=6,style=numeric,sorting=nyt]{biblatex}

\usepackage{hyperref}
\usepackage{microtype}

\newcommand{\De}[0]{\Delta}
\newcommand{\Ga}{\Gamma}
\newcommand{\F}{{\mathbb F}}
\newcommand{\C}{{\mathbb C}}
\newcommand{\R}{{\mathbb R}}
\newcommand{\N}{{\mathbb N}}

\newcommand{\delete}[1]{}

\newtheorem{Theorem}{Theorem}

\newtheorem{Proposition}[Theorem]{Proposition}
\newtheorem{Lemma}[Theorem]{Lemma}
\newtheorem{Problem}[Theorem]{Problem}
\newtheorem{innercustomthm}{Theorem}
\newenvironment{customthm}[1]
  {\renewcommand\theinnercustomthm{#1}\innercustomthm}
  {\endinnercustomthm}
 
\title{Classification by girth of three-dimensional algebraically defined monomial graphs over the real numbers}

\author[1]{Alex Kodess}
\author[2]{Brian G.\ Kronenthal}
\author[3]{Diego Manzano-Ruiz}
\author[4]{Ethan Noe}
\affil[1]{Farmingdale State College\\Farmingdale, NY 11735}
\affil[2,3,4]{Kutztown University of Pennsylvania\\Kutztown, PA 19530}
\begin{document}

\maketitle

\begin{abstract}
	For positive integers $s,t,u,v$, we define a bipartite graph $\Ga_\R(X^s Y^t,X^u Y^v)$ where each partite set is a copy of $\R^3$, and a vertex $(a_1,a_2,a_3)$ in the first partite set is adjacent to a vertex $[x_1,x_2,x_3]$ in the second partite set if and only if  
	\[
    	a_2 + x_2 = a_1^s x_1^t
    	\quad 
    	\text{and}
    	\quad
    	a_3+x_3=a_1^ux_1^v.
	\]  In this paper, we classify all such graphs according to girth.\\
	
\end{abstract}

\section{Introduction}\label{Section:Intro}

An algebraically defined graph in three dimensions $\Gamma_{\mathcal{R}}(f_2(X,Y),f_3(X,Y))$
is constructed using a ring $\mathcal{R}$ and bivariate functions $f_2,f_3\colon\mathcal{R}^2\to\mathcal{R}$.  These graphs are bipartite where each partite set is a copy of $\mathcal{R}^3$.  We label the vertices in the first partite set $(a_1,a_2,a_3)$ and in the second $[x_1,x_2,x_3]$.  In order for two vertices to be adjacent, denoted $(a_1,a_2,a_3)\sim[x_1,x_2,x_3]$, their coordinates must satisfy the equations $a_i+x_i=f_i(a_1,x_1)$ for $i=2,3$.
	
Dmytrenko, Lazebnik, and Williford \cite{dlw} studied the case where $\mathcal{R}$ is a finite field $\F_q$ of odd order and $f_2$ and $f_3$ are monomials (these graphs are aptly named monomial graphs).  They conjectured that all such monomial graphs of girth at least eight are isomorphic to $\Gamma_{\F_q}(XY,XY^2)$.  This work was expanded upon by Kronenthal \cite{k}, and the conjecture was ultimately proven by Hou, Lappano, and Lazebnik \cite{hll}.  In addition, Kronenthal and Lazebnik \cite{kl} and Kronenthal, Lazebnik, and Williford \cite{klw} studied families of polynomial graphs over algebraically closed fields of characteristic zero and applied some of their techniques to graphs over finite fields; these results were recently extended by Xu, Cheng, and Tang \cite{xu_tang}.  Moreover, Ganger, Golden, Kronenthal, and Lyons \cite{GGKL} studied a two-dimensional analogue over the real numbers. 
A number of questions related to connectivity, diameter, and isomorphisms of similarly constructed directed graphs, as well as a peculiar result on the number of roots of certain polynomials in finite fields, were considered in 
Kodess \cite{kod14}, Kodess and Lazebnik \cite{KL2015,KL2017}, 
Kodess, Lazebnik, Smith, and Sporre \cite{KLSS}, and  Coulter, De Winter, Kodess, and Lazebnik \cite{CDKL}.

In this paper, as in \cite{GGKL}, we study undirected graphs over the real numbers; but here, we examine the three-dimensional case.  Our main result is the following classification of all such monomial graphs  (the semicolons indicate logical conjunctions):

\begin{Theorem}\label{Thm:3DMonomials}
	Let $\Ga=\Ga_\R(X^s Y^t,X^u Y^v)$, where  $s,t,u,v\in\N$.  Then:
	\begin{enumerate}
		\item $\Ga$ has girth four if and only if at least one of $s,t$ is even and at least one of $u,v$ is even.
		\item $\Ga$ has girth six if and only if $\Ga$ is isomorphic to either $\Ga_\R(X^{2j+1} Y^{2k+1},X^{2m+1} Y^{2n})$ with $j\neq m$ or $n \leq k$, $\Ga_\R(X^{2j} Y^{2k},X^{2m+1} Y^{2n+1})$, or $\Ga_\R(X^{2j+1} Y^{2k+1},X^{2m+1} Y^{2n+1})$.
		
		In other words, $\Ga$ has girth six if and only if one of the following seven conditions is satisfied:
		\begin{enumerate}
			\item $s,t$ are even and $u,v$ are odd
			\item $s,t$ are odd and $u,v$ are even
			\item $s,t,u,v$ are odd
			\item $s$ is even; $t,u,v$ are odd; $t\neq v$ or $s< u$ 
			\item $t$ is even; $s,u,v$ are odd; $s\neq u$ or $t< v$ 
			\item $u$ is even; $s,t,v$ are odd; $t\neq v$ or $u< s$ 
			\item $v$ is even; $s,t,u$ are odd;  $s\neq u$ or $v<t$ 
		\end{enumerate}
		\item $\Ga$ has girth eight if and only if $\Ga$ is isomorphic to $\Ga_\R(X Y^{2k+1},X Y^{2n})$ with $n>k$.  
		
		In other words, $\Ga$ has girth eight if and only if one of the following four conditions is satisfied:
		\begin{enumerate}
			\item $s$ is even; $t,u,v$ are odd; $t=v$; $s> u$
			\item $t$ is even; $s,u,v$ are odd; $s=u$; $t> v$
			\item $u$ is even; $s,t,v$ are odd; $t=v$; $u> s$
			\item $v$ is even; $s,t,u$ are odd; $s=u$; $v> t$
		\end{enumerate}
	\end{enumerate}
\end{Theorem}

Before continuing, we make a few comments about Theorem \ref{Thm:3DMonomials}.  First, Theorem \ref{Thm:3DMonomials} is indeed a complete classification, i.e., every monomial graph in three dimensions is accounted for.  This is straightforward to check by considering the individual cases listed, but may be less obvious when examining the isomorphism classes presented in parts 2 and 3 until studied in concert with the isomorphisms from Lemmas \ref{lemma:graphisomorphisms} and \ref{Lemma:XisOne}.  Second, as $s$, $t$, $u$, and $v$ are positive integers, it immediately follows that $j$, $k$, $m$, and $n$ are non-negative integers and are strictly positive whenever they appear in an even exponent (e.g., $2j$ implies that $j$ must be a positive integer in that context, while $2j+1$ allows for the possibility $j=0$).

The study of algebraically defined graphs was originally motivated by several goals.  First, they are related to the construction of dense graphs of high girth, details of which can be found in Lazebnik, Sun, and Wang \cite{LSW} (Sections 4.4, 5.4, 6, and 7) and references therein.  The second motivation relates to incidence geometry, which we briefly discuss here.  In two dimensions, it is known (see Dmytrenko \cite{DmytrenkoThesis} and Lazebnik and Thomason \cite{lt}) that every graph $\Gamma_{\F_q}(f)$ with girth greater than four can be completed to a projective plane of order $q$ (although not all projective planes of order $q$ can be constructed in this way).  The three-dimensional analogue is motivated by the construction of generalized quadrangles because when $q$ is even, there exist monomial (and non-monomial) graphs $\Gamma_{\F_q}(f_2, f_3)$ that can be used to construct non-isomorphic generalized quadrangles of order $q$.  This motivated exploration of the $q$ odd case, as the existence of non-isomorphic girth eight algebraically defined graphs could lead to the construction of new generalized quadrangles.  For a detailed explanation of this construction, see \cite{kl} (concluding remarks), \cite{LSW} (Sections 4.2 and 4.3), and references therein.  This study was later generalized to algebraically defined graphs over fields of characteristic zero, including the field of complex numbers (as in \cite{kl,klw}) and the field of real numbers (as in this paper).

\section{Preliminary Tools \& Notation}\label{Section:Preliminaries}

We will begin by stating a necessary and sufficient condition for the existence of a $4$-cycle $(a, a_2,a_3)\sim [x, x_2,x_3]\sim(b, b_2,b_3)\sim[y, y_2, y_3]\sim(a, a_2, a_3)$ or a $6$-cycle $(a,a_2,a_3)\sim[x, x_2,x_3]\sim(b,b_2,b_3)\sim[y,y_2,y_3]\sim(c,c_2,c_3)\sim[z, z_2,z_3]\sim(a,a_2,a_3)$ in $\Gamma_\mathcal{R}(f_2,f_3)$.

\begin{Lemma}\label{lemma:delta}\textnormal{\cite{DmytrenkoThesis}}
	A $4$-cycle exists in $\Gamma_\mathcal{R}(f_2,f_3)$ if and only if there exist $a,b,x,y \in \mathcal{R}$ such that $a\ne b$, $x\ne y$, and for $i=2,3$,
	\begin{equation}f_i(a,x)-f_i(b,x)+f_i(b,y)-f_i(a,y)=0.\label{equation:delta2}
	\end{equation}
Similarly, a 6-cycle exists in $\Gamma_\mathcal{R}(f_2,f_3)$ if and only if there exist  distinct $a,b,c\in \mathcal{R}$ and distinct $x,y,z\in \mathcal{R}$ such that for $i=2,3$,
\begin{equation}
f_i(a,x)-f_i(b,x)+f_i(b,y)-f_i(c,y)+f_i(c,z)-f_i(a,z)=0.\label{equation:delta3}
\end{equation}

\end{Lemma}
\begin{proof}
This is straightforward from the 
definition of 
$\Ga = \Ga_{\cal R}(f_2,f_3)$: 
if vertices 
$(a,a_2,a_3)$, 
$[x,x_2,x_3]$, 
$(b,b_2,b_3)$, and 
$[y,y_2,y_3]$ 
are the consecutive vertices of a 
$4$-cycle in $\Ga$, then 
for $i = 2,3$,
\begin{equation}
\label{equation:4cycle}
\begin{aligned}
    a_i+ x_i &= f_i(a,x) \\ 
    b_i+ x_i &= f_i(b,x) \\ 
    b_i+ y_i &= f_i(b,y) \\
    a_i+ y_i &= f_i(a,y),
\end{aligned}
\end{equation}
and (\ref{equation:delta2}) follows. Conversely, if some four vertices satisfy equations (\ref{equation:4cycle}) above,  then they are the consecutive vertices of a $4$-cycle in $\Ga$, provided that the first coordinates of the vertices from the same partite sets are distinct (otherwise, if, for instance, $a = b$, then $a_2 = b_2$ and $a_3 = b_3$, and so the vertices $(a,a_2,a_3)$ and $(b,b_2,b_3)$ coincide.)

The argument for $6$-cycles is completely analogous and is therefore omitted.

\end{proof}

Since (\ref{equation:delta2}) and (\ref{equation:delta3}) appear repeatedly throughout this paper, we will introduce the following notation used, e.g., in \cite{DmytrenkoThesis,dlw,GGKL,kl,klw}:
\[\Delta_2(f_i)(a,b;x,y) 
:=
f_i(a,x)-f_i(b,x)+f_i(b,y)-f_i(a,y),\] and
\[
\Delta_3(f_i)(a,b,c;x,y,z)
:=
f_i(a,x)-f_i(b,x)+f_i(b,y)-f_i(c,y)+f_i(c,z)-f_i(a,z).
\]

Of particular interest, $\Delta_2(f_i)(a,b;x,y)$ and $\Delta_3(f_i)(a,b,c;x,y,z)$ depend only on the first coordinates of the vertices in the cycle.  Moreover, note that there will be many $4$-cycles with the same first coordinates as in a given $4$-cycle $(a,a_2,a_3)\sim[x,x_2,x_3]\sim(b,b_2,b_3)\sim[y,y_2,y_3]\sim(a,a_2,a_3)$, and we say that they are all of type $(a,b;x,y)$.  There will also be many $6$-cycles with the same first coordinates, so we say that they are all of type $(a,b,c;x,y,z)$.
For $k>3$, 
the definition of $\Delta_{2k}(f_i)(a_1,\ldots,a_k;x_1,\ldots,x_k)$ and 
notation for $2k$-cycle types are analogous to those of $k = 2$ and $k = 3$. We note that in cycles of length more than six, the first 
coordinates of the vertices in the same partite set need not be distinct.

We end this section with the following isomorphisms of the graph $\Gamma_\F(f_2,f_3)$, where $\F$ is a field; see, e.g., \cite{lt} (p.~3) or \cite{kl} (Proposition 2.2, p.~190) for proofs. First note that for a function $f=f(X,Y)$, we define $f^*$ as $f^*(X,Y) \coloneqq f(Y,X).$

\begin{Lemma}{\label{lemma:graphisomorphisms}}	
	Let $\F$ be a field and $f_2,f_3\in\F[X,Y]$. Then
	\begin{align}
	\Gamma_{\F}(f_2,f_3)&\cong\Gamma_{\F}(f_2^*,f_3^*),\label{equation:Iso1}\tag{$\mathcal{I}_1$}\\
	\Gamma_{\F}(f_2,f_3)&\cong\Gamma_{\F}(f_2,cf_3),\;\text{for all } c\in\F\backslash\{0\}, \label{equation:Iso2}\tag{$\mathcal{I}_2$}\\
	\Gamma_{\F}(f_2,f_3)&\cong\Gamma_{\F}(f_2,f_3+g+h),\;\text{for all}\; g\in\F[X] \;\text{and}\; h\in\F[Y]\label{equation:Iso3}\tag{$\mathcal{I}_3$}, \\
	 \Gamma_{\F}(f_2,f_3)&\cong\Gamma_{\F}(f_3,f_2)\label{equation:Iso4}\tag{$\mathcal{I}_4$}, \text{ and }\\
	 \Gamma_{\F}(f_2,f_3)&\cong\Gamma_{\F}(f_2,f_3+\delta f_2)\label{equation:Iso5}\tag{$\mathcal{I}_5$},\; \text{for any $\delta\in\F$}.
	\end{align}
\end{Lemma}

We will use (\ref{equation:Iso1}) to assume a given condition applies to $X$ instead of $Y$. By (\ref{equation:Iso4}), we can freely reverse the order of our monomials.  Moreover, by (\ref{equation:Iso2}), we assume without loss of generality throughout this paper that both $f_2$ and $f_3$ in $\Ga_{\R}(f_2,f_3)$ are of the form $X^iY^j$ for some $i,j\in\N$.  Finally, by virtue of isomorphisms (\ref{equation:Iso3}) and (\ref{equation:Iso5}), the results of Theorem \ref{Thm:3DMonomials} can be extended to certain families of graphs $\Gamma_{\F}(f_2,f_3)$, where $f_2$ and $f_3$ need not both be monomials.

\section{Proof of Theorem \ref{Thm:3DMonomials}}\label{Section:MainThm}

In this section, we prove Theorem \ref{Thm:3DMonomials} by breaking it up into a number of results.  The reader may note that while the statements of Theorem \ref{Thm:3DMonomials} are biconditionals, the results below are not worded as such.  However, Theorem \ref{Thm:3DMonomials} holds as stated because all graphs of the form $\Gamma_{\R}(X^sY^t,X^uY^v)$ are accounted for.  We begin by classifying all such graphs with girth four.

\begin{customthm}{\ref{Thm:3DMonomials} Part 1}
Let $\Ga=\Ga_{\R}(X^s Y^t,X^u Y^v)$ such that at least one of $s,t$ is even and at least one of $u,v$ is even.  Then $\Ga$ has girth four.
\end{customthm}

\begin{proof}
The graph $\Ga$ contains a 4-cycle of type $(1,-1;1,-1)$, as it is straightforward to check that (\ref{equation:delta2}) is satisfied.

\end{proof}

Note that in all remaining cases, either $s$ and $t$ are both odd or $u$ and $v$ are both odd.  We will now prove that all such graphs do not contain any 4-cycles.

\begin{Lemma}\label{Lemma:OddNoC4}
Let $\Ga=\Ga_{\R}(X^s Y^t,X^u Y^v)$ such that either $s$ and $t$ are both odd or $u$ and $v$ are both odd.  Then $\Ga$ has girth at least six.
\end{Lemma}

\begin{proof}
Suppose without loss of generality that $s$ and $t$ are both odd.  Then the graph $\Ga$ does not contain any 4-cycle $S$ of type $(a,b;x,y)$ because $\De_2(X^s Y^t)(S)=\left(a^s-b^s\right)\left(x^t-y^t\right)=0$ has no real solutions satisfying $a\neq b$ and $x\neq y$.
\end{proof}
For proofs in the remainder of this paper, Lemma \ref{Lemma:OddNoC4} will allow us to skip over the proof that each graph is 4-cycle free, and instead focus on the presence (or absence) of larger cycles.  We next examine a family of girth six graphs.

\begin{customthm}{\ref{Thm:3DMonomials} Part 2abc}
Let $\Ga=\Ga_{\R}(X^s Y^t,X^u Y^v)$ such that $s$ and $t$ have the same parity, $u$ and $v$ have the same parity, and $s$, $t$, $u$, and $v$ are not all even.  Then $\Ga$ has girth six.
\end{customthm}

\begin{proof}
    The graph $\Ga$ contains a 6-cycle of type $(1,0,-1;-1,1,0)$, as it is straightforward to check that (\ref{equation:delta3}) is satisfied.  By Lemma \ref{Lemma:OddNoC4}, the result follows.  
\end{proof}

We briefly pause to comment that Lemma \ref{Lemma:OddNoC4} and Theorem \ref{Thm:3DMonomials} Part 2abc do not hold over $\C$.  For example, the graph $\Gamma_{\R}(X^3Y^3,X^6Y^6)$ has girth six, but $\Gamma_{\C}(X^3Y^3,X^6Y^6)$ has girth four because it contains a 4-cycle of type $\left(1,-\frac{1}{2}+\frac{\sqrt{3}}{2}i\ ;1,-1\right)$.  

We will next prove the remaining cases of Theorem \ref{Thm:3DMonomials} Part 2.

\begin{customthm}{\ref{Thm:3DMonomials} Part 2defg}
If $\Ga = \Ga_\R(X^{2j+1} Y^{2k+1},X^{2m+1} Y^{2n})$ such that $j\neq m$ or $n\leq k$, then $\Ga$ has girth six.
\end{customthm}

This theorem is an immediate consequence of Propositions \ref{Prop-Thm2defg-1}, \ref{Prop-Thm2defg-2}, and \ref{Prop-Thm2defg-3}, which follow below.  Before stating them, we recall that for a positive real number $a$, arbitrary integer $b$ and  positive odd integer $c$, the exponent $(-a)^{\frac{b}{c}}$ is unambiguously and uniquely defined as 
$(-1)^b \sqrt[c]{a^b}$.

\begin{Proposition} \label{Prop-Thm2defg-1}
	If $n\leq k$, then the graph $\Ga=\Ga_\R(X^{2j+1} Y^{2k+1},X^{2m+1} Y^{2n})$ has girth six.
\end{Proposition}
\begin{proof}
We show that $\Ga$ contains a 6-cycle $S$ of type $(0,1,-1; 1,y,z)$ for some $y,z\in\R$ with $y$, $z$, and $1$ all distinct, which occurs if and only if 
$\De_3\left(X^{2j+1}Y^{2k+1}\right)(S) = 
\De_3\left(X^{2m+1}Y^{2n}\right)(S) = 0$. 
These equations yield 
$z=\left(2y^{2k+1}-1\right)^{\frac{1}{2k+1}}$, 
and therefore 
$D(y):=
-1+2y^{2n}-\left(2y^{2k+1}-1\right)^{\frac{2n}{2k+1}}=0
$. 
The condition $n \le k$ ensures that $D$ has a 
root $y \in (-\infty,0)$ since $D(0) = -2$ and ${\displaystyle \lim_{y\to-\infty} D(y) = \infty}$.

It is now easy to see that $y$, $z=\left(2y^{2k+1}-1\right)^{\frac{1}{2k+1}}$, and $1$ are distinct, and so $\Ga$ contains a $6$-cycle $S$ of this type.

\end{proof}

From now on we can assume that $n > k$ and note that $j \neq m$ is equivalent to $s \neq u$.

\begin{Proposition}\label{Prop-Thm2defg-2}
	If $n>k$ and $m<j$, then the graph $\Ga=\Ga_\R(X^{2j+1} Y^{2k+1},X^{2m+1} Y^{2n})$ has girth six.
\end{Proposition}

\begin{proof}
We show that $\Ga$ contains a 6-cycle $S$ of type $(0,1,c; 0,1,z)$ for some $c,z\in\R$ with $c$, $z$, and $1$ all distinct, which occurs if and only if 
$\De_3\left(X^{2j+1}Y^{2k+1}\right)(S) = 
\De_3\left(X^{2m+1}Y^{2n}\right)(S) = 0$. 
These equations yield 
$z=\left(\frac{c^{2j+1}-1}{c^{2j+1}}\right)^{\frac{1}{2k+1}}$, 
and therefore 
$D(c):=1-c^{2m+1}+c^{2m+1}\left(\frac{c^{2j+1}-1}{c^{2j+1}}\right)^{\frac{2n}{2k+1}} = 0$. 
The condition $n > k$ ensures that $D$ has a root $c \in (-\infty,-1)$ since $D(-1) < 0$ and, by elementary calculus,  ${\displaystyle \lim_{c\to-\infty} D(c) = 1}$ as $j > m$.

It is now easy to see that $c$, $z=\left(\frac{c^{2j+1}-1}{c^{2j+1}}\right)^{\frac{1}{2k+1}}$, and $1$ are distinct, and so $\Ga$ contains a $6$-cycle $S$ of this type.

\end{proof}

\begin{Proposition}\label{Prop-Thm2defg-3}
	If $n>k$ and $m>j$, then the graph $\Ga=\Ga_\R(X^{2j+1} Y^{2k+1},X^{2m+1} Y^{2n})$ has girth six.
\end{Proposition}
\begin{proof}
We show that $\Ga$ contains a 6-cycle $S$ of type $(1,0,-3; x,1,z)$ for some $x,z\in\R$ with $x$, $z$, and $1$ all distinct, which occurs if and only if $\De_3\left(X^{2j+1}Y^{2k+1}\right)(S) = \De_3\left(X^{2m+1}Y^{2n}\right)(S) = 0$. 
These equations yield 
$z=\left(\frac{x^{2k+1}+3^{2j+1}}{3^{2j+1}+1}\right)^{\frac{1}{2k+1}}$, 
and therefore 
$D(x):=	x^{2n}+3^{2m+1}-\left(3^{2m+1}+1\right)\left(\frac{x^{2k+1}+3^{2j+1}}{3^{2j+1}+1}\right)^{\frac{2n}{2k+1}}=0$. The derivative $D'$ is continuous at $x=1$, and we note that $D(1) = 0$, and also
$D'(1)=2n\left(1-\frac{3^{2m+1}+1}{3^{2j+1}+1}\right)<0$ as $m>j$. Hence, there exist $x_1,x_2\in\R$ with $x_1 < 1 < x_2$ such that $D(x_2) < 0 < D(x_1)$. We also note that 
$\frac{3^{2m+1}+1}{(3^{2j+1}+1)^{\frac{2n}{2k+1}}} \neq 1 $ since otherwise 
$3^{2m+1}+1 = 
	\left(
	\frac{(3^{2j}+1)^n}{(3^{2m+1}+1)^k}
	\right)^2$ 
	is a square, which is not possible: if for some integer $t$, we have $3^{2m+1} + 1 = t^2$, then $t-1$ and $t+1$ are both powers of $3$, so $t=2$,  contradicting $m>j\ge 0$. This ensures that either 
 ${\displaystyle\lim_{x\to\pm\infty}D(x) = \infty}$ or ${\displaystyle \lim_{x\to\pm\infty}D(x) = -\infty}$, and so $D$ has a root $x$ either in $(-\infty,x_1)$ or in $(x_2,\infty)$.

It is now easy to see that $x$, $z=\left(\frac{x^{2k+1}+3^{2j+1}}{3^{2j+1}+1}\right)^{\frac{1}{2k+1}}$, and $1$ are distinct, and so $\Ga$ contains a $6$-cycle $S$ of this type.

\end{proof}

The following lemma builds on Lemma \ref{lemma:graphisomorphisms} by introducing an additional isomorphism of real monomial graphs.  This isomorphism will be used in the proof of Theorem \ref{Thm:3DMonomials} Part 3.

\begin{Lemma}\label{Lemma:XisOne}
	The monomial graphs $\Ga_\R(X^{2m+1} Y^t,X^{2m+1} Y^v)$ and $\Ga_\R(X Y^t,X Y^v)$ are isomorphic, where $m$ is a non-negative integer and $t,v$ are positive integers.
\end{Lemma}

\begin{proof}
	This is a result of the following isomorphism:
	
	\begin{align*}
	\Ga_\R(X^{2m+1} Y^t,X^{2m+1} Y^v)&\to\Ga_\R(X Y^t,X Y^v)\\
	(x_1,x_2,x_3)&\mapsto (x_1^{2m+1},x_2,x_3)\\
	[y_1,y_2,y_3]&\mapsto[y_1,y_2,y_3]
	\end{align*}
\end{proof}

\begin{customthm}{{\ref{Thm:3DMonomials} Part 3}}
If $n>k$, then the graph $\Ga_\R(X^{2m+1} Y^{2k+1},X^{2m+1} Y^{2n})$ has girth eight.
\end{customthm}

\begin{proof}
By the previous lemma, $\Ga_\R(X^{2m+1} Y^{2k+1},X^{2m+1} Y^{2n})\cong \Ga_\R(X Y^{2k+1},X Y^{2n})$.  Therefore, we only need to prove that $\Ga=\Ga_\R(X Y^{2k+1},X Y^{2n})$ has girth eight.

Suppose that for some distinct $a,b,c\in\R$ and some  distinct $x,y,z\in\R$ the graph $\Ga$ contains a 6-cycle $S$ of type $(a,b,c;x,y,z)$, and so
	\begin{align*}
	\De_3(XY^{2k+1})(S)&=x^{2k+1}(a-b)+y^{2k+1}(b-c)+x^{2k+1}(c-a)=0\\
	\De_3(XY^{2n})(S)&=x^{2n}(a-b)+y^{2n}(b-c)+x^{2n}(c-a)=0.
	\end{align*}
We thus obtain a homogeneous linear system in variables $a-b$, $b-c$, and $c-a$, with matrix
\[
A=
\begin{bmatrix}
	1&1&1\\
	x^{2k+1}&y^{2k+1}&z^{2k+1}\\
	x^{2n}&y^{2n}&z^{2n}\\
\end{bmatrix}.
\]
We will show that  
\[
\det A=
\left(
y^{2k+1}-x^{2k+1}\right)\left(z^{2n}-x^{2n}\right)-\left(z^{2k+1}-x^{2k+1}\right)\left(y^{2n}-x^{2n}\right)
\]	
is nonzero for any distinct $x$, $y$, and $z$, implying a desired contradiction $a=b=c$. We have $\det A \neq 0$ if and only if 
\begin{equation}\label{eqn:monotone}
	\frac{z^{2n}-x^{2n}}{z^{2k+1}-x^{2k+1}}
	\neq
	\frac{y^{2n}-x^{2n}}{y^{2k+1}-x^{2k+1}}
	\,
	\text{ for all distinct } 
	x,y,z\in\R.
\end{equation}
We prove (\ref{eqn:monotone}) by verifying that, for any $x\in\R$, the function 
	\[
	h_x(t):=\frac{t^{2n}-x^{2n}}{t^{2k+1}-x^{2k+1}},\quad 
	t\in\R\setminus\{x\},
	\]
is strictly increasing on its entire domain. This is true for $h_0(t)=t^{2n-2k-1}$, and for any $x\neq 0$, consider the derivative 
\[
	h_x'(t)
	=\frac{t^{2k}\left((2n-2k-1)t^{2n}-2nx^{2k+1}t^{2n-2k-1}+(2k+1)x^{2n}\right)}{\left(t^{2k+1}-x^{2k+1}\right)^2}.\label{Equation:hprime}
\] 
Define 
$H_x(t):=(2n-2k-1)t^{2n}-2nx^{2k+1}t^{2n-2k-1}+(2k+1)x^{2n}$, $t\in\R$. It is easy to see that $H_x(t) > 0$ for all $t\in \R \setminus \{x\}$ because
$H_x$ is continuous on $\R$, 
its only critical points are $t = 0$ and $t = x$ with 
$H_x(0) = (2k+1)x^{2n} > H_x(x) = 0$, and 
 ${\displaystyle \lim_{t\to\pm\infty}H_x(t) = \infty}$. It follows that $h'_x(0) = 0$ and $h'_x(t) > 0$ if $t\in\R\setminus\{0,x\}$. We also note that ${\displaystyle \lim_{t\to x^-}h_x(t) = \lim_{t\to x^+}h_x(t)\in\R}$. Thus $h_x(t)$ is strictly increasing for any $x\in\R$.
 
 It is straightforward to check that $\Ga$ has an $8$-cycle of type $(1,0,-1,0; 1,-1,1,-1)$, and so $\Ga$ has girth eight.
\end{proof}

\section{Concluding Remarks}\label{Section:Conclusion}
The primary goal of study in \cite{hll} was to ascertain that when ${\mathcal R}$ is a finitehttps://www.overleaf.com/project/5db074df0a3a250001cea841 field $\F_q$ of odd order $q$, the only (up to isomorphism) girth eight graph $\Ga_{{\mathcal R}}(f,g)$, where $f$ and $g$ are  monomials in ${\mathcal R}[X,Y]$, is $\Ga_{{\mathcal R}}(XY,XY^2)$. 
A similar assertion was proven in \cite{kl} and \cite{klw}: whenever ${\mathcal R}$ is an algebraically closed field of characteristic zero, the only (up to isomorphism) graph $\Ga_{{\mathcal R}}(X^k Y^m,g)$ of girth at least eight, where $k,m\in\N$ and $g\in{\mathcal R}[X,Y]$, is $\Ga_{{\mathcal R}}(X Y,XY^2)$. It was proven \cite{xu_tang} that 
given any polynomials 
$f\in{\mathbb F}_q[X]$, $g\in{\mathbb F}_q[Y]$, and $h\in{\mathbb F}_q[X,Y]$, 
there exists a positive integer $M$ depending on the degrees of $f$, $g$, and $h$, such that 
any graph $\Ga_{{\mathcal R}}(fg,h)$ with ${\mathcal R} = \F_{q^M}$ of girth at least eight is isomorphic to 
$\Ga_{{\mathcal R}}(XY,XY^2)$; it was also proven that when $\mathcal{R}$ is  any algebraically closed field of characteristic zero, the only graph $\Ga_{{\mathcal R}}( f(X) g(Y), h(X,Y))$ (up to isomorphism) of girth at least eight is 
$\Ga_{{\mathcal R}}( XY, X Y^2)$. 
Finally, there are no graphs of girth eight or more in the two-dimensional real case; see \cite{GGKL}. This discussion motivates the following question about girth eight graphs in the three-dimensional real case.

\begin{Problem}
Sort by isomorphism the graphs $\Ga_{{\mathbb R}}(XY^{2k+1},XY^{2n})$, $n > k$.
\end{Problem}

Of particular interest is the question whether all such graphs are isomorphic to $\Ga_{{\mathbb R}}(XY,XY^2)$; this stands in contrast to the situation  ${\cal R} = \F_q$ with $q$ odd (as in \cite{hll}) and ${\cal R} = \C$ (as in \cite{kl,klw}), where it is known that all monomial graphs of girth eight are isomorphic to $\Ga_{\cal R}(XY, XY^2)$. We do not for instance know whether $\Ga_{{\mathbb R}}(XY,XY^2)$ is isomorphic to $\Ga_{{\mathbb R}}(XY,XY^4)$.

We note that the automorphism group of $\Ga_{{\mathbb R}}(XY,XY^2)$ acts transitively on each of the partite sets $\{(p_1,p_2,p_3)\}$ and $\{[l_1,l_2,l_3]\}$; also $\Ga_{{\mathbb R}}(XY,XY^2)$ is edge-transitive. These statements are easily established in the wake of the following automorphisms of $\Ga_{{\mathbb R}}(XY,XY^2)$:

    \begin{align}
        \begin{split}
            (p_1,p_2,p_3) &\mapsto 
            (p_1+a,
            p_2,
            p_3), \\
            [l_1,l_2,l_3] &\mapsto
            [l_1,
            l_2+al_1,
            l_3+al_1^2
            ],\quad a \in\R.
            \label{aut:type1}
        \end{split}        
        \tag{${\mathcal A}_1$}\\[2ex]
        \begin{split}
        (p_1,p_2,p_3) &\mapsto 
        (p_1,
        p_2+bp_1,
        p_3+2bp_2+b^2p_1), \\
        [l_1,l_2,l_3] &\mapsto
        [l_1+b,
        l_2,
        l_3+2bl_2
        ],\quad b \in\R.
        \label{aut:type2}
        \end{split}        
        \tag{${\mathcal A}_2$}\\[2ex]        
        \begin{split}
        (p_1,p_2,p_3) &\mapsto 
        (p_1,p_2-c,p_3-d),\\
        [l_1,l_2,l_3] &\mapsto
        [l_1,l_2+c,l_3+d], 
        \quad c,d\in\R.
        \label{aut:type3}
        \end{split}
        \tag{${\mathcal A}_3$}
    \end{align}

We do not know whether $\Ga_{{\mathbb R}}(XY,XY^4)$ has any of the three transitivity properties mentioned above, and while $\Ga_{{\mathbb R}}(XY,XY^4)$ certainly has automorphisms of type (\ref{aut:type3}), it is not clear how to construct automorphisms of type (\ref{aut:type1}) or (\ref{aut:type2}) for it.

\bigskip

\noindent \textbf{Acknowledgments:} 
The authors are thankful to the anonymous referees whose thoughtful comments improved
the paper, 
to Tony Wong for a useful suggestion regarding the proof of Theorem 1 Part 3,  
and to Felix Lazebnik for fruitful discussions related to the isomorphism question discussed in Section \ref{Section:Conclusion}.
The work of the first author was supported by Farmingdale State College, and the work of the fourth author was supported by the Kutztown University Bringing Experiences About Research in Summer Program.  

\printbibliography
\typeout{get arXiv to do 4 passes: Label(s) may have changed. Rerun}
\end{document}